\documentclass[11pt]{amsart}
\usepackage{amsmath,amssymb, amsthm}

\newtheorem{theorem}{Theorem}[section]
\newtheorem{corollary}[theorem]{Corollary}
\newtheorem{lemma}[theorem]{Lemma}
\newtheorem{proposition}[theorem]{Proposition}

\theoremstyle{definition}

\newcommand*{\be}{\begin{equation}}
\newcommand*{\ee}{\end{equation}}

\newcommand{\nn}{\nonumber}
\providecommand{\abs}[1]{\lvert#1\rvert}

\newcommand*{\fl}[1]{\lfloor{#1}\rfloor}

\def\phat{\widehat{p}}
\def\qhat{\widehat{q}}
\def\Qhat{\widehat{Q}}

\def\ind{\mathbf{1}}
\DeclareMathOperator{\Var}{Var}
\DeclareMathOperator{\Cov}{Cov}

\def\bN{\mathbb{N}}
\def\bZ{\mathbb{Z}}

\def\la{\lambda}

\def\omhat{\widehat\omega}
\def\rhohat{\widehat\rho}
\def\omtil{\widetilde\omega}

\def\pcount{I}
\def\range{R}

\def\ratea{p}
\def\rateb{q}

\def\flux{H}  

\def\Qlb{m_Q}  
\def\lbla{a}  
\def\lblb{b} 
\def\mom{m} 

\newcommand*{\Ev}{\mathbf E}
\newcommand*{\Vv}{{\text{\bf Var}}} 
\newcommand*{\Cv}{{\text{\bf Cov}}}
\newcommand*{\Pv}{\mathbf P}  
\newcommand*{\hop}{\bigskip\noindent}

\newcommand*{\om}{\omega}

\newcommand*{\Zb}{\mathbb Z}

\numberwithin{equation}{section}
\allowdisplaybreaks[1]

\begin{document}

\title[Fluctuations for ASEP]{Fluctuation bounds for the 
asymmetric simple exclusion process}  
\author[M.~Bal\'azs]{M\'arton Bal\'azs}
\address{M\'arton Bal\'azs,
Budapest University of Technology and Economics, 
Institute of Mathematics, 1 Egry J\'ozsef u., H \'ep.\ V.7., 
Budapest, Hungary}
\email{balazs@math.bme.hu}
\urladdr{http://www.math.bme.hu/~balazs}
\thanks{M. Bal\'azs was partially supported by the Hungarian Scientific Research Fund (OTKA) grants K60708, TS49835, F67729, and the Bolyai Scholarship of the Hungarian Academy of Sciences.}
\author[T.~Sepp\"al\"ainen]{Timo Sepp\"al\"ainen}
\address{Timo Sepp\"al\"ainen\\ University of Wisconsin-Madison\\ 
Mathematics Department\\ Van Vleck Hall\\ 480 Lincoln Dr.\\  
Madison WI 53706-1388\\ USA.}
\email{seppalai@math.wisc.edu}
\urladdr{http://www.math.wisc.edu/~seppalai}
\thanks{T.\ Sepp\"al\"ainen was partially supported by National Science Foundation
grant DMS-0701091.} 
\keywords{exclusion process, particle current, second class particle,
coupling, variance bound, moment bound, fluctuations}
\subjclass[2000]{60K35} 
\date{\today}
\begin{abstract}
We give a partly new proof of the fluctuation bounds for the second
class particle and current in the stationary 
asymmetric simple exclusion process. One novelty is a coupling
that preserves the ordering of second class particles in two
systems that are themselves ordered coordinatewise. 
\end{abstract}
\maketitle


\section{Introduction}
The  asymmetric
simple exclusion process (ASEP) is a Markov process that
describes the motion of 
  particles on the one-dimensional 
integer lattice $\Zb$. Each particle executes 
a continuous-time nearest-neighbor random walk on $\bZ$
with jump rate $p$ to the right and $q$ to the left.  
Particles interact through the exclusion rule
which  means that 
 at most one particle is allowed at each site.  
Any attempt to jump onto an already occupied site is
ignored.  The asymmetric case is $p\ne q$.
We  assume  $0\le q< p\le 1$ 
and also $p+q=1$.
 
In this paper we give a partially new proof of the 
fluctuation bounds for ASEP first proved in 
\cite{se2/3} and \cite{varj2nd}. 
For a more thorough explanation of the context and
related work we refer the reader to these earlier
papers.  The present paper has two
novelties.

\medskip

(i) Both the earlier and the present 
approach are based on identities that
link together the current, space-time covariances, and
the second class particle. In \cite{varj2nd} we proved these
identities with  martingale techniques 
and generator computations.  Section \ref{sc:covar} of the 
present paper shows that  these identities are 
 consequences of not much more than 
simple  counting of particles.  These arguments are 
elementary and  should work very generally.   Without much 
effort we extend the identities to more
general exclusion dynamics that allow 
 bounded jumps and rates that
depend on the local configuration. However, product invariant 
distributions remains a key assumption. 

\medskip

(ii) In \cite{se2/3} we used sophisticated couplings introduced in
\cite{fks} and some complicated estimation to bound the 
positions of certain second class particles.  
In Section \ref{sc:newcoup} of the present paper we
introduce a coupling for ASEP that keeps the second class particle 
of a denser system behind the second class particle of a 
comparison system.  Since the macroscopic speed of a second
class particle is $\flux'(\rho)$, the slope of the flux at 
the given density $\rho$,  concavity of the flux suggests that
second class particles travel slower in a denser system. 
Since this is what the coupling achieves, we think of 
this coupling as a form of {\sl microscopic concavity.}

\medskip

Fix two parameters $0\le q<p\le 1$ such that $p+q=1$. 
We run quickly through the fundamentals 
of  $(p,q)$-ASEP.
We refer the reader
to standard references for further details 
 \cite{ips, stochi}.

\medskip

{\bf Definition and graphical construction.}
ASEP represents the motion of particles on the
integer lattice $\bZ$.  The particles are subject to
an exclusion rule which means that each site of
$\bZ$ can contain at most one particle.
   The state of the system at time $t$ is
a configuration $\eta(t)=(\eta_i(t))_{i\in\bZ}\in\{0,1\}^\bZ$
of zeroes and ones.  The value $\eta_i(t)=1$ means that 
 site $i$ is occupied by a particle at time $t$, while the
value  $\eta_i(t)=0$ means that 
 site $i$ is vacant at time $t$. 

 The motion of the particles is controlled by
independent Poisson processes ({\sl Poisson clocks})  
$\{N^{i\to i+1}, N^{i\to i-1}: i\in\bZ\}$.  
Each Poisson clock  $N^{i\to i+1}$ has rate $p$
 and each  $N^{i\to i-1}$ has rate $q$. If
$t$ is  a jump time for $N^{i\to i+1}$ and if
$(\eta_i(t-), \eta_{i+1}(t-))=(1,0)$ then at time
$t$ the particle from site $i$ moves to site $i+1$
and the new values are $(\eta_i(t), \eta_{i+1}(t))=(0,1)$.
Similarly if $t$ is  a jump time for $N^{i\to i-1}$ a 
particle is moved from $i$ to $i-1$ at time $t$,
provided the configuration at time $t-$ permits this
move.  If the jump prompted by a Poisson clock is not
permitted by the state of the system, this jump attempt
 is simply
ignored and the particles resume waiting for the next
prompt coming from the Poisson clocks.   

This construction of the process is known as the {\sl graphical
construction} or the {\sl Harris construction}. 
When the initial state is a fixed configuration $\eta$,  
$P^\eta$ denotes the distribution of the process. 
 
We write $\eta$, $\omega$, etc for elements of the state space
$\{0,1\}^\bZ$, but also for the entire process so that 
$\eta$-process
 stands for $\{\eta_i(t): i\in\bZ, 0\le t<\infty\}$. 
The configuration $\delta_i$ is the state that has a single
particle at position $i$ but otherwise the lattice is
vacant.  

\medskip

{\bf Invariant distributions.} 
A basic fact is that i.i.d.~Bernoulli distributions 
$\{\nu^\rho\}_{\rho\in[0,1]}$ are 
extremal invariant distributions for ASEP. 
For each density value $\rho\in[0,1]$, $\nu^\rho$ is the 
probability measure on $\{0,1\}^\bZ$ 
 under which the occupation variables 
$\{\eta_i\}$ are i.i.d.\ with common mean 
$\int \eta_i\,d\nu^\rho=\rho$.  When the process
$\eta$ is stationary with time-marginal $\nu^\rho$, we write
$P^\rho$ for the probability distribution of the entire
process.  
The {\sl stationary density-$\rho$ process} means the ASEP
$\eta$ that is stationary in time and has marginal
distribution $\eta(t)$ $\sim$ $\nu_\rho$. 

\medskip

{\bf Basic coupling and second class particles.}
The {\sl basic coupling} of two exclusion processes
$\eta$ and $\om$ means that they obey a common set
of Poisson clocks $\{N^{i\to i+1},N^{i\to i-1}\}$. 
Suppose the two processes $\eta$ and $\eta^+$ satisfy 
$\eta^+(0)=\eta(0)+\delta_{Q(0)}$ at time zero, for 
some position $Q(0)\in\bZ$. This means that $\eta^+_i(0)=\eta_i(0)$
for all $i\ne Q(0)$, $\eta^+_{Q(0)}(0)=1$ and $\eta_{Q(0)}(0)=0$.
  Then throughout the evolution in the basic coupling
there is a single discrepancy between $\eta(t)$ and  $\eta^+(t)$
at some position $Q(t)$:  $\eta^+(t)=\eta(t)+\delta_{Q(t)}$.
From the perspective of $\eta(t)$, $Q(t)$ is called a 
second class particle. By the same token,
from the perspective of $\eta^+(t)$ $Q(t)$ is a 
second class {\sl anti}particle.
  In particular, we shall call the 
pair $(\eta, Q)$ {\sl a $(p,q)$-ASEP with a second class particle}. 

We write a boldface $\Pv$ for the probability measure
when more than one process are coupled together. 
In particular, $\Pv^\rho$ represents the
situation where the initial occupation 
variables  $\eta_i(0)=\eta^+_i(0)$
are i.i.d.~mean-$\rho$ Bernoulli for $i\ne 0$, and the
second class particle $Q$ starts at $Q(0)=0$.   

More generally, if two processes $\eta$ and $\om$ are in basic
coupling and $\om(0)\ge\eta(0)$ (by which we mean coordinatewise
ordering $\om_i(0)\ge\eta_i(0)$ for all $i$) then the ordering
$\om(t)\ge\eta(t)$ holds for all $0\le t<\infty$.  The effect of
the basic coupling is to give priority to the $\eta$ particles
over the $\om-\eta$ particles. Consequently we can think
of the $\om$-process as consisting of  first class particles 
(the $\eta$ particles) and second class particles 
(the $\om-\eta$ particles).

\medskip

{\bf Current.} 
For $x\in\bZ$ and $t>0$,
$J_x(t)$ stands for the net left-to-right
particle current across the straight-line space-time path 
from $(1/2,0)$ to $(x+1/2,t)$.  More precisely, 
$J_x(t)=J_x(t)^+-J_x(t)^-$ where 
$J_x(t)^+$ is the number of particles that lie
in $(-\infty,0]$ at time $0$ but lie in 
$[x+1,\infty)$ at time $t$, while 
$J_x(t)^-$ is the number of particles that lie
in $[1,\infty)$ at time $0$ and  in 
 $(-\infty,x]$ at time $t$. When more than one process
($\om$, $\eta$, etc)
is considered in a coupling, the currents of the processes
are denoted by $J^\om_x(t)$, $J^\eta_x(t)$, etc. 

\medskip

{\bf Flux and characteristic speed.}  The average net rate at which
particles in the stationary  $(p,q)$-ASEP at density $\rho$
 move across a fixed edge $(i,i+1)$ is the {\sl flux}
\be
\flux(\rho)=(p-q)\rho(1-\rho).
\label{def:flux}\ee
The characteristic speed at density $\rho$ is 
\be
V^\rho=\flux'(\rho)=(p-q)(1-2\rho).
\label{def:Vrho}\ee
In the stationary process the expected  current
is 
\be
E^\rho[J_x(t)]=t\flux(\rho)-x\rho 
\label{eq:EJ}\ee
as can be seen by noting that particles that crossed the
edge $(0,1)$  either also crossed $(x,x+1)$ and 
contributed to $J_x(t)$ or did not. 
Another important and less obvious expectation is 
\be
\Ev^\rho[Q(t)]=t\flux'(\rho)
\label{eq:EQH'a}\ee 
for the second class particle that starts at the origin. 
We derive \eqref{eq:EQH'a} in Section \ref{sc:covar}.

We can now state the main result, the moment bounds on the
second class particle.  

\begin{theorem}
Consider  $(p,q)$-ASEP with rates such that 
$0\le q=1-p<p\le 1$.  

{\rm (Upper bound)} 
For densities $\rho\in(0,1)$  there
exists a  constant $0<C_1(\rho)<\infty$ such that 
 for $1\le \mom<3$ and  $t\ge 1$
\be
\Ev^\rho\bigl\{\,\lvert Q(t)-V^\rho t\rvert^\mom\,\bigr\}  
 \;\le\; \frac{C_1(\rho)}{3-\mom} t^{2\mom/3}. 
\label{eq:mainubd}\ee
$C_1(\rho)$ is a continuous function of $\rho$.

{\rm (Lower bound)} 
For densities $\rho\in(0,1)$ there
exist constants $0<t_0(\rho), C_2(\rho)<\infty$ 
such that for  $t\ge t_0(\rho)$ and all $\mom\ge 1$ 
\be
\Ev^\rho\bigl\{\,\lvert Q(t)-V^\rho t\rvert^\mom\,\bigr\}  
 \;\ge\; C_2(\rho)t^{2\mom/3}. 
\label{eq:mainlbd}\ee
$C_2(\rho)$ and $t_0(\rho)$ are continuous
functions of $\rho\in(0,1)$.   
\label{th:2classmom}\end{theorem}

The constants  $C_1, t_0, C_2$ depend also on $p$
but since $p$ is regarded as fixed we omit this dependence
from the notation.  The boundary dependence on $\rho$ 
 is as follows: as   $\rho\to\{0,1\}$,    
$C_1(\rho)\to\infty$, 
 $C_2(\rho)\to 0$  and 
$t_0(\rho)\to\infty$. 
The convergence $C_2(\rho)\to 0$ as $\rho\to\{0,1\}$
is natural because for $\rho\in\{0,1\}$
$Q(t)$ is a random walk.  
The divergences of  $C_1(\rho)$ and $t_0(\rho)$ 
appear to be artifacts of our proof. 

For the third moment and beyond
 our upper bound argument gives these bounds: for $0<\rho<1$
\be
\Ev^\rho\bigl\{\,\lvert Q(t)-V^\rho t\rvert^\mom\,\bigr\}  \le 
\begin{cases}
B(\rho,3) t^2\log t &\text{for $\mom=3$ and $t\ge e$,}\\
{B(\rho,m)} t^{\mom-1} &\text{for $3<\mom<\infty$ and $t\ge 1$.}
\end{cases}\label{eq:addub}\ee 
The constant $B(\rho,m)$ is finite and  
  continuous in $\rho$ for each fixed $3\le \mom<\infty$. 

A key identity proved in the next section states that 
\[
\Var^\rho[J_{\fl{V^\rho t}}(t)]
=\rho(1-\rho)\Ev^\rho\lvert Q(t)-\fl{V^\rho t}\rvert.\]
From this and Theorem \ref{th:2classmom}
 we get the order of the variance
of the current as seen by an observer traveling at the 
characteristic speed $V^\rho$: for large enough $t$, 
\[
C_1(\rho)t^{2/3} \le \Var^\rho[J_{\fl{V^\rho t}}(t)] 
\le C_2(\rho) t^{2/3}.
\]

In the case of TASEP ({\sl totally} asymmetric simple
exclusion process with $p=1=1-q$)  P.~L.~Ferrari and Spohn
\cite{ferspohn}
proved a distributional limit for the current 
$J_{\fl{V^\rho t}}(t)$.  At the time of this writing
this precise distributional limit has not yet been proved for 
ASEP. 
If the observer chooses a speed $V\ne V^\rho$ then
$\Var^\rho[J_{\fl{Vt}}(t)]$ is of order $t$ and in fact
$J_{\fl{Vt}}(t)$ satisfies a central limit theorem
that  records Gaussian fluctuations of the
initial particle configuration \cite{se}.

\medskip

{\bf Organization of the paper.}  Section 
\ref{sc:covar} proves the identities that connect the 
 current and the second class particle. 
In Section \ref{sc:newcoup} we develop the coupling that
keeps two second class particles in different densities 
ordered.  Section \ref{sc:RW} contains a tail bound
on a biased random walk in an inhomogeneous environment 
that is needed for the last section. 
Section \ref{sc:2momproof} proves Theorem \ref{th:2classmom}. 

\medskip

{\bf Further notation.}  $\bZ_+=\{0,1,2,\dotsc\}$
and $\bN=\{1,2,3,\dotsc\}$.  
  Constants denoted by $C$ or $C_i$ ($i=1,2,3,\dotsc$) 
can change from line to line. Centering of a random variable
is denoted by   $\widetilde X=X-EX$.  Shift on
$\{0,1\}^\bZ$ is denoted by $(\theta_i\om)_j=\om_{i+j}$.

\section{Covariance identities}
\label{sc:covar}

In this section we can consider more general exclusion
processes with bounded jumps and rates that depend
on the configuration around the jump location. 

Let 
$\range$ be a constant  that bounds the range of
admissible  jumps. 
 For $1\le k\le \range$ let $p_k$ and $q_k$ be functions
of particle configurations that depend only
on coordinates $(\om_i: -\range\le i\le \range)$ 
and satisfy 
$0\le \ratea_k(\om), \rateb_k(\om)\le 1$. 
The rule for the evolution is that for all $i\in\bZ$ and 
all $k\in\{1,\dotsc,\range\}$, independently of everything
else,  the  exchange $\om\mapsto\om^{i,i+k}$ happens
 at rate 
\[ \ratea_k(\theta_i\om)\om_i(1-\om_{i+k})
+ \rateb_k(\theta_{i+k}\om)\om_{i+k}(1-\om_{i}). \]
The configuration $\om^{i,i+k}$ is the result of 
exchanging the contents of sites $i$ and $i+k$:
\[
\om^{i,i+k}_j=\begin{cases}  \om_{i+k}, &\text{if $j=i$}\\
                         \om_{i}, &\text{if $j=i+k$}\\
                 \om_{j}, &\text{if $j\notin\{i,i+k\}.$}
\end{cases} \]
The rule of evolution can be restated as follows: 
 whenever possible a particle jumps 
from $i$ to $i+k$ at rate $\ratea_k(\theta_i\om)$ and 
from $i$ to $i-k$ at rate $\rateb_k(\theta_i\om)$. 

\medskip

{\bf Key assumption.} Bernoulli distributions 
$\{\nu^\rho\}_{\rho\in[0,1]}$ are invariant for the 
process.  As throughout the paper, $P^\rho$, $E^\rho$,
$\Var^\rho$ and $\Cov^\rho$ 
refer to the density-$\rho$ invariant process.

\medskip

We state formulas that tie together the current, space-time
covariance (also called ``two-point function''),
and the second class particle.  We begin with the well-known
formula that connects the two-point function with the 
second class particle. 

Let $\Pv^\rho$ be the probability distribution of two
coupled processes $\om\le \om^+$ that start with identical 
Bernoulli-$\rho$ occupation variables 
$\om_i(0)=\om^+_i(0)$ at $i\ne 0$ and a single discrepancy
$Q$  that starts at the origin.  In other words
$Q(0)=0$ and $\om^+(t)=\om(t)+\delta_{Q(t)}$. 
Then for $0<\rho<1$ 
\be
\Cov^\rho[\om_j(t),\,\om_0(0)]
=\rho(1-\rho) \Pv^\rho\{Q(t)=j\}.
\label{eq:CovQ}\ee
Note that the two sides of the identity come from different
processes: the left-hand side is a covariance in a stationary
process, while the right-hand side is in terms of processes
perturbed at the origin at time $0$. 
For the sake of completeness we give below a proof of 
\eqref{eq:CovQ}. 

\begin{theorem} For any density $0\le \rho\le 1$, 
$z\in\bZ$ and $t>0$ we have these formulas:  
\be
\Var^\rho[J_z(t)]=\sum_{j\in\bZ} \abs{j-z} \Cov^\rho[\om_j(t),\om_0(0)]
=\rho(1-\rho)\Ev^\rho\lvert Q(t)-z\rvert. 
\label{goal1}\ee
and 
\be
\frac{d}{d\rho} E^\rho[J_z(t)]
= \Ev^\rho [Q(t)]-z. 
\label{eq:ddtEJ}\ee
The terms in the series in {\rm\eqref{goal1}} decay
exponentially in $\abs{j}$, uniformly over $\rho$. 
All members 
of these identities  are continuous functions of 
$\rho\in[0,1]$. At $\rho=0$ and $\rho=1$  the left-hand side of 
 {\rm\eqref{eq:ddtEJ}} is a one-sided derivative. 
\label{th:formulas}
\end{theorem} 

Some comments and consequences follow. 

At $\rho=0$ the second class
particle  $Q(t)$ is a random walk that takes
jumps of size  $k$ at rate $p_k(\delta_0)$
and  jumps of size  $-k$ at rate $q_k(\delta_0)$.
At the other extreme,  $\rho=1$, 
  $Q(t)$ is a random walk that takes
jumps of size  $k$ at rate $q_k(\mathbf{1}-\delta_{-k})$
and  jumps of size  $-k$ at rate $p_k(\mathbf{1}-\delta_{k})$.
Here $\mathbf{1}$ denotes the configuration $\om\equiv 1$. 
In particular,  in $(p,q)$-ASEP, 
as $\rho$ goes from $0$ to $1$, $Q$ 
interpolates between 
nearest-neighbor random walks with rates $(p,q)$ and $(q,p)$. 

The equilibrium current past the origin satisfies 
\be
E^\rho[J_0(t)]=t\flux(\rho)
\label{eq:EJH}\ee
where the flux $H(\rho)$ is the expected rate of particle
motion across any fixed edge in the stationary density-$\rho$
process: 
\[
H(\rho)=\sum_{k=1}^\range kE^\rho(p_k-q_k).
\]
 Combining   \eqref{eq:ddtEJ} for $z=0$ and 
\eqref{eq:EJH} 
gives the useful identity 
\be
\Ev^\rho[Q(t)]=t\flux'(\rho).
\label{eq:EQH'}\ee 

We turn to the proofs, beginning with \eqref{eq:CovQ}.  

\begin{proof}[Proof of equation \eqref{eq:CovQ}]
This is a straight-forward calculation. 
\begin{align*}
&\Cov^\rho[\om_j(t),\om_0(0)]
=E^\rho[\om_j(t)\om_0(0)]-\rho^2
=\rho E^\rho[\om_j(t)\,\vert\,\om_0(0)=1]-\rho^2\\
&=\rho\Bigl( E^\rho[\om_j(t)\,\vert\,\om_0(0)=1]
-\rho E^\rho[\om_j(t)\,\vert\,\om_0(0)=1]
-(1-\rho) E^\rho[\om_j(t)\,\vert\,\om_0(0)=0]\Bigr)\\
&=\rho(1-\rho)\Bigl( E^\rho[\om_j(t)\,\vert\,\om_0(0)=1]
- E^\rho[\om_j(t)\,\vert\,\om_0(0)=0]\Bigr)\\
&=\rho(1-\rho)\Bigl( \Pv^\rho[\om^+_j(t)=1]
- \Pv^\rho[\om_j(t)=1] \Bigr)
=\rho(1-\rho) \Pv^\rho[\om^+_j(t)=1\,,\, \om_j(t)=0] \\
&=\rho(1-\rho) \Pv^\rho[ Q(t)=j].
\qedhere\end{align*}
\end{proof} 

\bigskip

The remainder of this section proves 
 Theorem \ref{th:formulas}.  
The second equality in \eqref{goal1} comes from 
\eqref{eq:CovQ}. 

Let $\om$ be a stationary exclusion process satisfying the 
assumptions made in this section, with i.i.d.~Bernoulli($\rho$) distributed
occupations $\{\om_i(t)\}$ at any fixed time $t$. 
To approximate the infinite system with
finite systems, for each $N\in\bN$ 
let process $\om^N$ have initial configuration
\be
\om^N_i(0)=\om_i(0) \ind_{\{-N\le i\le N\}}.
\label{Ninitial}\ee 
We assume that all these processes are coupled through
 a Harris-type construction, with jump attempts prompted by 
Poisson clocks, with appropriate  rates, 
attached to directed edges $(i,j)$ for $\abs{j-i}\le \range$.
Let $J^N_z(t)$ denote the current in process $\om^N$.

Let $z(0)=0$, $z(t)=z$, and introduce the counting variables
\be
\pcount^N_+(t)=\sum_{n> z(t)}\om^N_n(t)\,,\quad
\pcount^N_-(t)=\sum_{n\le z(t)}\om^N_n(t).
\label{defpcount}\ee
Then the current can be expressed as 
\[
J_z^N(t)= \pcount^N_+(t)-\pcount^N_+(0)=\pcount^N_-(0)-\pcount^N_-(t),
\]
and its variance as 
\begin{align*}
\Var J_z^N(t) &=\Cov\bigl(\pcount^N_+(t)-\pcount^N_+(0), 
\,\pcount^N_-(0)-\pcount^N_-(t)\bigr)\\
&=\Cov\bigl(\pcount^N_+(t), \pcount^N_-(0))  
+ \Cov(\pcount^N_+(0),\pcount^N_-(t))\\
&\qquad -\Cov(\pcount^N_+(0),\pcount^N_-(0)) 
-\Cov(\pcount^N_+(t),\pcount^N_-(t))\\
&=\sum_{k\le 0,\,m>z}  \Cov[\om^N_m(t),\om^N_k(0)]
+ \sum_{k\le z,\, m>0}  \Cov[\om^N_k(t),\om^N_m(0)]\\
&\qquad\qquad\qquad -\Cov(\pcount^N_+(0),\pcount^N_-(0)) 
-\Cov(\pcount^N_+(t),\pcount^N_-(t)).
\end{align*}

Independence gives
\[\Cov(\pcount_+^N(0),\pcount_-^N(0))=0\] and the identity
above simplifies
to
\be
\begin{split}
\Var J_z^N(t) &=
\sum_{k\le 0,\,m>z}  \Cov[\om^N_m(t),\om^N_k(0)]\\
&\qquad + \sum_{k\le z,\,m>0}  \Cov[\om^N_k(t),\om^N_m(0)]
-\Cov(\pcount_+^N(t),\pcount_-^N(t)).
\end{split}
\label{simplified}\ee 
We  show that identity \eqref{simplified} converges 
to identity \eqref{goal1}. 
 
To take advantage of the decaying correlations that result
from the bounded jump range, we introduce another family
of auxiliary processes. 
Let $\eta^{a,b}$ and  $\eta^{N,a,b}$ be
 exclusion processes defined from initial conditions
\be
\eta^{a,b}_i(0)=\om_i(0) \ind_{\{a< i<b\}}
\quad\text{and}\quad
\eta^{N,a,b}_i(0)=\om^N_i(0) \ind_{\{a< i<b\}}
\label{etaabinit}\ee
and with a ``reduced'' Poisson construction:
 all Poisson jumps that involve
any site outside $(a,b)$ are deleted.
 The point of this definition
is that for disjoint intervals $(a,b)$ and $(u,v)$, processes
 $\eta^{a,b}$  and   $\eta^{u,v}$ are independent
because they do not share initial occupation variables or
Poisson clocks. 

We write $\Pv$ for the probability measure under which all
these coupled processes $\{\om,\om^N,$ $\eta^{a.b},\eta^{N,a,b}\}$
live.   The next lemma is valid for completely general 
initial occupations  $\{\om_i(0)\}$. 

\begin{lemma} Let $\{\om_i(0)\}$ be an arbitrary deterministic
or random initial configuration, 
and define initial configurations
 $\om^N(0)$, $\eta^{a,b}(0)$ and 
$\eta^{N,a,b}(0)$ 
 by {\rm\eqref{Ninitial}} and {\rm\eqref{etaabinit}}. 
For a fixed $0<t<\infty$ there is a constant 
$C(t)<\infty$ 
such that, for all indices $N$, $m$, $k$,
 and all times $s\in[0,t]$ 
\begin{align}
\Ev\lvert \om^N_k(s)- \eta^{N,k-m,k+m}_k(s)\rvert &\le e^{-C(t)m}
\label{aux11}\\
\text{and}\quad
\Ev\lvert \om_k(s)- \eta^{k-m,k+m}_k(s)\rvert &\le e^{-C(t)m}.
\label{aux11.5}\end{align}
\label{etaauxlm} 
\end{lemma}

\begin{proof}
We give the argument for a particular $N$ and  $k$. 
But the argument and the resulting bound have
no dependence on   $N$, $k$ or the  initial configurations. 
The only issue is  the disconnectedness of 
the graph created by the Poisson clocks in the Harris construction
during time interval $[0,t]$. 

For any given integer $i>0$, there is a fixed positive probability
that no site in the intervals 
 $[k-i-\range+1, k-i]$  and $[k+i,k+i+\range-1]$
 is involved in any Poisson jump
during time interval $[0,t]$.   Hence the probability that 
this event fails at each integer $i=3\range \ell$ 
for $1\le \ell\le n$ 
is exponentially small in $n$. 
(For distinct $\ell$ the clocks involved are independent.) 
 
If this event succeeds for $i=3\range\ell$, then 
$\om^N_k(s)=\eta^{N,k-m,k+m}_k(s)$ for all  
$m\ge 3\range\ell+\range$ and $s\in[0,t]$
 by the following reasoning.
This $m$ is large enough to  ensure that the initial configurations  
$\om^N(0)$ and $\eta^{N,k-m,k+m}(0)$ agree  on
$(k-i-\range, k+i+\range)$.   
 The absence of Poisson
jumps has frozen the configurations in 
$(k-i-\range, k-i]\cup[k+i, k+i+\range)$ up to time $t$.
Also, up to time $t$
no Poisson jump connects 
the process inside the interval 
$(k-i, k+i)$ to the outside. 
The rates 
 of the jumps
inside the interval $(k-i, k+i)$ are determined by the configuration
 in the larger interval
$(k-i-\range, k+i+\range)$. 
Consequently inside $(k-i,k+i)$, up to time $t$, the processes 
$\om^N$ and $\eta^{N,k-m,k+m}$ execute the same
moves. 
\end{proof}

Now combine the previous lemma with independent initial
occupations. 

\begin{corollary} Let $\{\om_i(0)\}$ be   i.i.d.~Bernoulli($\rho$) and 
 $\om^N(0)$ defined   by 
 {\rm\eqref{Ninitial}}.
Fix $0<t<\infty$. 
 There exists a constant $C(t)$ 
such that  for all indices $N,i,k$, all 
densities $0\le \rho\le 1$ and times $s\in[0,t]$, 
\be
\bigl\lvert\,\Ev\omtil^N_i(s)\omtil^N_k(t)\,\bigr\rvert
\le 4e^{-C(t)\abs{k-i}}
\label{aux12}\ee
and 
\be
\bigl\lvert\,\Cov^\rho[\om_i(s),\om_k(t)]\,\bigr\rvert
\;=\;\bigl\lvert\,\Ev\omtil_i(s)\omtil_k(t)\,\bigr\rvert
\le 4e^{-C(t)\abs{k-i}}.
\label{aux14}\ee

 For any function $g$ on $\bZ$ 
that grows at most polynomially, 
\be
[0,1]\ni \rho\mapsto \sum_k g(k)\Cov^\rho[\om_k(t),\om_0(0)]
\label{covregf}\ee
is a continuous function.  
\label{auxcor1}\end{corollary}

\begin{proof} Let $m=\abs{i-k}$. 
Variables $\eta^{N,i-m/2,i+m/2}_i(s)$ and 
$\eta^{N,k-m/2, k+m/2}_{i+m}(t)$
are independent and so their covariance is zero.
\eqref{aux12} follows from \eqref{aux11}. Similarly for 
\eqref{aux14}. 

For fixed $k,m$ the function
\be
 \rho\mapsto  \Cov^\rho[\eta^{k-m,k+m}_k(t),\om_0(0)]
\label{aux14.1}\ee
is continuous because the expectation
 depends on $\rho$ through  the finitely many
variables $\{\eta^{k-m,k+m}_i(0): k-m<i<k+m, \, \om_0(0)\}$.
These can be coupled simultaneously for all values
of $\rho$ through i.i.d.\ uniform $U_i\sim$ Unif$(0,1)$
by writing 
\be
\eta^{k-m,k+m,\rho}_i(0)=\ind\{0<U_i<\rho\}.  
\label{eq:Ucoup}\ee
Continuity in \eqref{aux14.1}  follows by dominated convergence. 
Estimate \eqref{aux11.5} then shows the continuity
of the individual terms in the series in \eqref{covregf}. 
The continuity of the whole series 
follows from the $\rho$-uniformity of the tail bounds in \eqref{aux14}.
\end{proof}

Once $N>\abs{a}\vee\abs{b}$,  $\eta^{N,a,b}(t)=\eta^{a,b}(t)$.
From this and Lemma \ref{etaauxlm}  one concludes that 
\be
\om^N_i(t)\to \om_i(t) \quad\text{as $N\to\infty$, a.s. and in 
$L^2$, for any $i$.}
\label{aux13}\ee

From these estimations we derive this lemma:

\begin{lemma} The right-hand side of {\rm\eqref{simplified}}
converges as $N\to\infty$ to the middle member
 of {\rm\eqref{goal1}}.
\label{auxlemma1}\end{lemma}

\begin{proof}
We begin by showing the convergence of the last term of
\eqref{simplified}: 
\be
\Cv(\pcount_+^N(t),\pcount_-^N(t))\to 0.
\label{aux16}\ee
This follows from writing, for a fixed $K>0$, 
\begin{align*} 
&\Cv(\pcount_+^N(t),\pcount_-^N(t))
= \sum_{k>z}\sum_{j\le z} \Cv\bigl(\om^N_k(t),\,\om^N_j(t)\bigr)\\
&= 
\sum_{0<k-z\le K}\sum_{-K\le j-z\le 0} 
\Cv\bigl(\om^N_k(t),\,\om^N_j(t)\bigr)\\
&\qquad +\sum_{k>z}\sum_{j<z-K} \Ev\omtil^N_k(t)\omtil^N_j(t)
+\sum_{k>z+K}\sum_{-K\le j-z\le 0} \Ev\omtil^N_k(t)\omtil^N_j(t)\\
&= 
\sum_{0<k-z\le K}\sum_{-K\le j-z\le 0} 
\Cv\bigl(\om^N_k(t),\,\om^N_j(t)\bigr)\\
&\qquad \qquad +O\Bigl(\; \sum_{k>z, j<z-K} e^{-C(k-j)} +
\sum_{k>z+K, -K\le j-z\le 0} e^{-C(k-j)}\Bigr)\\
&= 
\sum_{0<k-z\le K}\sum_{-K\le j-z\le 0} \Cv\bigl(\om^N_k(t),\,\om^N_j(t)\bigr)
+O( e^{-CK}).
\end{align*} 
For a fixed $K$ as $N\to\infty$ 
\[
\sum_{0<k-z\le K}\sum_{-K\le j-z\le 0} \Cv\bigl(\om^N_k(t),\,\om^N_j(t)\bigr)
\to
\sum_{0<k-z\le K}\sum_{-K\le j-z\le 0} \Cv\bigl(\om_k(t),\,\om_j(t)\bigr)
=0\]
where the vanishing is due to the assumption of i.i.d.~variables
at time $t$.  Letting $K\nearrow\infty$ finishes the proof of
\eqref{aux16}.

A similar argument, approximating with the variables 
$\eta^{N,\,k-\abs{k}/2\,,\,k+\abs{k}/2}_{{k}}(t)$ 
in the tails of the series and taking
the $N\to\infty$ limit in a finite sum, proves 
\begin{align}
&\sum_{k\le 0,\,m>z}  \Cv[\om^N_m(t),\om^N_k(0)]
+ \sum_{k\le z,\, m>0}  \Cv[\om^N_k(t),\om^N_m(0)] \label{covline7}\\
&\ \longrightarrow\ 
\sum_{k\le 0,\, m>z}  \Cv[\om_m(t),\om_k(0)]
+ \sum_{k\le z,\, m>0}  \Cv[\om_k(t),\om_m(0)]\nn\\
&=\sum_{n\in\bZ} \abs{n-z} \Cv[\om_n(t),\om_0(0)]. \nn
\end{align}
The last equality used shift-invariance. 
This proves Lemma \eqref{auxlemma1}. 
\end{proof} 

To complete the proof of \eqref{goal1}
it only remains to show that the 
left-hand side of {\rm\eqref{simplified}}
converges as $N\to\infty$ to the  first member of {\rm\eqref{goal1}}.
The same line of reasoning works again because 
up to an exponentially small probability in $m$, the current
is not altered by removal of all jumps that involve sites 
outside $(-m,m)$. 

Identity \eqref{goal1} is now proved. The
statement about the uniform convergence of the series and
the continuity in $\rho$ come from Corollary \ref{auxcor1}.

\bigskip

We turn to the proof of \eqref{eq:ddtEJ}. 

\begin{lemma} For densities $0<\la<\rho<1$, 
currents in stationary processes satisfy 
\be
E^\rho[J_z(t)]-E^\la[J_z(t)] = \int_\la^\rho \frac1{\theta(1-\theta)}
\sum_{j\in\bZ} (j-z) \Cov^\theta[\om_j(t),\,\om_0(0)]\,d\theta.
\label{eq:Jdiff}\ee
\label{Jdifflm1}\end{lemma}

\begin{proof} Let again process with superscript $N$ be the one whose
 initial configuration is \eqref{Ninitial} with i.i.d.\ occupations
on $[-N,N]$ and vacant sites elsewhere, and $\Pv$ the measure
for all the coupled processes. 
 Let $I^N=\sum_i\om^N_i(t)$ be the number of particles in the
process $\om^N$. $I^N$ is a Binomial($2N+1,\rho$) random variable. 
For $0<\rho<1$
\begin{align}
\frac{d}{d\rho} \Ev[J_z^N(t)] &=
\frac{d}{d\rho} \sum_{m=0}^{2N+1} \binom{2N+1}{m} 
\rho^m(1-\rho)^{2N+1-m} E[J^N(t)\vert I^N=m] \nn\\
&=
 \sum_{m=0}^{2N+1} P(I^N=m)
\Bigl(\frac{m}{\rho} -\frac{2N+1-m}{1-\rho}\Bigr) E[J^N(t)\vert I^N=m] \nn\\
&=\frac1{\rho(1-\rho)} \Ev\Bigl[ J_z^N(t)\cdot
\sum_i (\om^N_i(0)-\rho)\Bigr]\nn\\
\intertext{[recall definitions \eqref{defpcount}]}
&=\frac1{\rho(1-\rho)} \Cv\bigl[ \pcount^N_+(t)-\pcount^N_+(0)\,,\,
\pcount^N_-(0)+\pcount^N_+(0)\bigr]\nn\\
&=\frac1{\rho(1-\rho)}\Bigl( \Cv\bigl[ \pcount^N_+(t)\,,\,
\pcount^N_-(0)\bigr]+\Cv\bigl[ \pcount^N_+(t)-\pcount^N_+(0)\,,\,
\pcount^N_+(0)\bigr]\Bigr). 
 \label{diffJline4} 
\end{align}
The last equality used $\Cv[\pcount^N_+(0),\pcount^N_-(0)]=0$
that comes from the i.i.d.\ distribution of initial occupations. 
The first covariance on line \eqref{diffJline4} write directly
as 
\[
\Cv\bigl[ \pcount^N_+(t)\,,\,
\pcount^N_-(0)\bigr]
= 
\sum_{k\le 0,\,m>z}  \Cv[\om^N_m(t),\om^N_k(0)]. 
\]
The second 
 covariance on line \eqref{diffJline4} write as 
\begin{align*} 
&\Cv\bigl[ \pcount^N_+(t)-\pcount^N_+(0)\,,\,
\pcount^N_+(0)\bigr]
= \Cv\bigl[ \pcount^N_-(0)-\pcount^N_-(t)\,,\,
\pcount^N_+(0)\bigr]\\
&\qquad =\;-\Cv\bigl[ \pcount^N_-(t)\,,\,
\pcount^N_+(0)\bigr]
=\ - \sum_{k\le z,\, m>0}  \Cv[\om^N_k(t),\om^N_m(0)]. 
\end{align*}
Inserting these back on line  \eqref{diffJline4} gives
\[
\frac{d}{d\rho} \Ev[J_z^N(t)] =  
\frac1{\rho(1-\rho)}\Bigl( \;
\sum_{k\le 0,\,m>z}  \Cv[\om^N_m(t),\om^N_k(0)] \ -\ 
\sum_{k\le z,\, m>0}  \Cv[\om^N_k(t),\om^N_m(0)] \;
\Bigr).
\]
Thus compared to line \eqref{covline7} we have the difference
instead of the sum.  Integrate over the density $\rho$ and
take $N\to\infty$ as was taken on the line following \eqref{covline7}.
This proves the proposition.   
\end{proof}

This lemma together with \eqref{eq:CovQ} 
 gives 
\be
E^\rho[J_z(t)]-E^\la[J_z(t)] = \int_\la^\rho 
(\Ev^\theta[Q(t)]-z)\,d\theta 
\label{eq:Jdiff1}\ee
for $0<\lambda<\rho<1$.
Couplings show the continuity of these expectations:
\[
E^\la[J_z(t)]\to E^\rho[J_z(t)]
\quad\text{and}\quad
\Ev^\lambda[Q(t)]\to \Ev^\rho[Q(t)]
\quad\text{as $\lambda\to\rho$ in $[0,1]$.}
\]
Precisely speaking, if the initial configurations
are  coupled as indicated in \eqref{eq:Ucoup},  the integrands
converge a.s.\ to the corresponding integrands.   Then
bounds in terms of Poisson processes give uniform integrability 
that makes the expectations converge.

Thus the right-hand side of \eqref{eq:Jdiff1} can be
differentiated in $\rho$ and 
identity \eqref{eq:ddtEJ} for $0<\rho<1$  follows.  

For the one-sided derivatives, consider the case $\rho=0$. 
If we take $\lambda\searrow 0$ in 
\eqref{eq:Jdiff1} then
\[
E^\rho[J_z(t)]= \int_0^\rho 
(\Ev^\theta[Q(t)]-z)\,d\theta.  
\]
Continuity of the integrand now allows us to differentiate
from the right at $\rho=0$.  Similar argument for the left derivative
at $\rho=1$ completes the proof of identity 
\eqref{eq:ddtEJ}. 
The continuity of the right-hand side of \eqref{eq:ddtEJ}
was also argued along the way. 

\section{A coupling for microscopic concavity}
\label{sc:newcoup} 

As observed in \eqref{eq:EQH'} the speed of the second class particle 
in a density-$\rho$ ASEP is $\flux'(\rho)$.  
Thus by the concavity of the flux $\flux$ 
 a defect travels slower  in a denser system (recall that 
we assume $p>q$ throughout).  However,
in ASEP the basic coupling does not respect this, except in the
totally asymmetric ($p=1, q=0$) case.  To see this,
consider two pairs of
processes $(\om^+,\om)$ and $(\eta^+,\eta)$ such that both
pairs have one discrepancy: $\om^+(t)=\om(t)+\delta_{Q^\om(t)}$
and $\eta^+(t)=\eta(t)+\delta_{Q^\eta(t)}$.  Assume that 
$\om(t)\ge \eta(t)$. 
In basic coupling  the jump from state
\[
\begin{bmatrix}  \om^+_i&\om^+_{i+1}\\
\om_i&\om_{i+1}\\ \eta^+_i&\eta^+_{i+1}\\ 
\eta_i&\eta_{i+1}\end{bmatrix} 
=\begin{bmatrix} 1&1\\0&1\\1&0\\0&0\end{bmatrix}
\quad\text{to state}\quad
\begin{bmatrix} 1&1\\1&0\\1&0\\0&0\end{bmatrix}
\]
happens at rate $q$ and results in $Q^\om=i+1>i=Q^\eta$.

In this section we construct  a different coupling that combines 
the basic coupling with auxiliary clocks for second
class particles.  The idea is to think of a single ``special''
second class particle as performing a random walk on the 
process of $\om-\eta$ second class particles.   This  coupling 
preserves the expected ordering of the special second class particles, 
 hence it can be regarded as a form 
of {\sl microscopic concavity}.

This theorem  summarizes the outcome. 

\begin{theorem} 
Assume given two  initial configurations
 $\{\zeta_i(0)\}$ and $\{\xi_i(0)\}$ and two
not necessarily distinct  positions 
$Q^\zeta(0)$ and $Q^\xi(0)$  on $\bZ$.
Suppose the coordinatewise ordering  $\zeta(0)\ge\xi(0)$ holds,
$Q^\zeta(0)\le Q^\xi(0)$, 
and $\zeta_i(0)=\xi_i(0)+1$ for $i\in\{Q^\zeta(0), Q^\xi(0)\}$.
Define the configuration $\zeta^-(0)=\zeta(0)-\delta_{Q^\zeta(0)}$. 

Then there exists a coupling of processes 
$(\zeta^-(t), Q^\zeta(t), \xi(t), Q^\xi(t))_{t\ge 0}$ with initial 
state $(\zeta^-(0), Q^\zeta(0), \xi(0), Q^\xi(0))$ given in 
the previous paragraph, such that both pairs 
$(\zeta^-, Q^\zeta)$ and $(\xi, Q^\xi)$ are $(p,q)$-ASEP's
with a second class particle, and $Q^\zeta(t)\le Q^\xi(t)$
for all $t\ge 0$. 
\label{th:newcoupling}
\end{theorem} 

To begin the construction, 
put two exclusion processes
$\zeta$ and $\xi$ in basic coupling, obeying 
Poisson clocks $\{N^{i\to i\pm 1}\}$. 
They are  ordered so that $\zeta\ge\xi$. 
The $\zeta-\xi$ second class particles are labeled
in increasing order 
$\dotsm<X_{m-1}(t)<X_{m}(t)<X_{m+1}(t)<\dotsm$. 
We assume  there is at least
one such second class particle, but
beyond that we make no assumption about their
number.  Thus there is 
 some
finite or infinite subinterval $I\subseteq\bZ$ of indices 
such that 
  the positions of the  $\zeta-\xi$ second class particles
are given by $\{X_{m}(t): m\in I\}$. 

We introduce two dynamically evolving labels 
$\lbla(t)$, $\lblb(t)$ $\in I$ in such a manner that 
$X_{\lbla(t)}(t)$ is the position of a second class 
antiparticle
in the $\zeta$-process, 
$X_{\lblb(t)}(t)$ is the position of a second class particle
in the $\xi$-process, and the ordering 
\be
X_{\lbla(t)}(t)\le X_{\lblb(t)}(t)
\label{eq:2classorder} \ee
 is preserved by the 
dynamics. 

The labels $\lbla(t)$, $\lblb(t)$ are allowed to jump from
$m$ to $m\pm 1$ 
only when particle $X_{m\pm 1}$ is adjacent  
 to $X_m$. The labels  do not take
jump commands from the
Poisson clocks $\{N^{i\to i\pm 1}\}$
that govern $(\xi,\zeta)$.  Instead, the directed
edges $(i,i+1)$ and $(i,i-1)$ are given another
collection of independent Poisson clocks so that the
following 
jump rates are realized. 

(i) If $\lbla=\lblb$ and $X_{\lbla+1}=X_{\lbla}+1$ 
then
\[
 (\lbla,\lblb)\ \text{jumps to }\ \begin{cases}
 (\lbla,\lblb+1) &\text{with rate $p-q$}\\ 
(\lbla+1,\lblb+1)&\text{with rate $q$.}
\end{cases} \]

(ii) If $\lbla=\lblb$ and $X_{\lbla-1}=X_{\lbla}-1$ 
then
\[
 (\lbla,\lblb)\ \text{jumps to }\ \begin{cases}
 (\lbla-1,\lblb) &\text{with rate $p-q$}\\ 
(\lbla-1,\lblb-1)&\text{with rate $q$.}
\end{cases} \]

(iii) If $\lbla\ne\lblb$ then $\lbla$ and $\lblb$ jump
independently with these rates:
\begin{align*}
 \lbla\ \text{jumps to }\ &\begin{cases}
 \lbla+1 &\text{with rate $q$ if $X_{\lbla+1}=X_{\lbla}+1$}\\ 
 \lbla-1 &\text{with rate $p$ if $X_{\lbla-1}=X_{\lbla}-1$;}
\end{cases} \\
 \lblb\ \text{jumps to }\ &\begin{cases}
 \lblb+1 &\text{with rate $p$ if $X_{\lblb+1}=X_{\lblb}+1$}\\ 
 \lblb-1 &\text{with rate $q$ if $X_{\lblb-1}=X_{\lblb}-1$.}
\end{cases} 
\end{align*}

Let us emphasize that the pair process $(\xi,\zeta)$ is still
governed by the old clocks  $\{N^{i\to i\pm 1}\}$
in the basic coupling. The new clocks that realize
rules (i)--(iii) are not observed except when sites
$\{i,i+1\}$ are both occupied by $X$-particles and 
at least one of $X_{\lbla}$ or $X_{\lblb}$ lies in
$\{i,i+1\}$.    

First note that if initially $\lbla(0)\le\lblb(0)$ then
 jumps (i)--(iii)  preserve the inequality
$\lbla(t)\le \lblb(t)$ which gives \eqref{eq:2classorder}.  
(Since the jumps in point (iii) happen
independently, there cannot be two simultaneous jumps. 
So it is not possible for $\lbla$ and $\lblb$ to
cross each other with a $(\lbla,\lblb)\to(\lbla+1,\lblb-1)$
move.)  

Define processes  $\zeta^-(t)=\zeta(t)-\delta_{X_{\lbla(t)}(t)}$
and  $\xi^+(t)=\xi(t)+\delta_{X_{\lblb(t)}(t)}$.
In other words, to produce $\zeta^-$ remove particle 
$X_{\lbla}$ from $\zeta$, and  to produce $\xi^+$ add particle 
$X_{\lblb}$ to $\xi$.  The second key point is that,
even though these new processes are no longer defined
by the standard graphical construction, distributionwise
 they are still  ASEP's with second class particles. 
We argue this point for  $(\zeta^-, X_{\lbla})$ and leave the argument 
for $(\xi, X_{\lblb})$ to the reader. 

\begin{lemma}
The pair $(\zeta^-, X_{\lbla})$ 
is a $(p,q)$-ASEP with a second class particle. 
\label{lm:zeta-asep}\end{lemma}

\begin{proof} We check that the jump rates for the 
process  $(\zeta^-, X_{\lbla})$, produced by the combined
effect of the basic coupling with 
 clocks  $\{N^{i\to i\pm 1}\}$ and the new clocks, 
are  the same jump rates that result from defining
an (ASEP, second class particle) pair in terms of the
graphical construction as explained in the Introduction.  

To have notation for the
possible jumps, let $0$ denote an empty site, $1$ 
a $\zeta^-$-particle, and $2$ particle $X_{\lbla}$. 
Consider a fixed pair $(i,i+1)$ of sites and write
$xy$ with $x,y\in\{0,1,2,\}$ for the contents of sites
$(i,i+1)$ before and after the jump.  Then here 
are the possible moves across
the edge $\{i,i+1\}$, and the rates that these moves
would have in the basic coupling. 
\begin{align*}
\text{Type 1} \qquad &10\longrightarrow 01 \ \text{ with rate $p$}\\
                  &01\longrightarrow 10 \ \text{ with rate $q$}\\[4pt]
\text{Type 2} \qquad &20\longrightarrow 02 \ \text{ with rate $p$}\\
                  &02\longrightarrow 20 \ \text{ with rate $q$}\\[4pt]
\text{Type 3} \qquad &12\longrightarrow 21 \ \text{ with rate $p$}\\
                  &21\longrightarrow 12 \ \text{ with rate $q$}\\
\end{align*}
Our task is to check that the construction of 
$(\zeta^-,X_{\lbla})$ actually
realizes these rates.  

Jumps of types 1 and 2 are prompted by
the clocks  $\{N^{i\to i\pm 1}\}$
of the graphical construction of $(\xi,\zeta)$, and hence have the correct
rates listed above.  

Jumps of type 3 occur in two distinct
ways. 

(Type 3.1) First there can be a $\xi$-particle next to $X_{\lbla}$, and then
the rates shown above are again realized by 
the clocks  $\{N^{i\to i\pm 1}\}$ because in
the basic coupling the $\xi$-particles have priority 
over the $X$-particles.  

(Type 3.2) The other 
alternative is 
that both sites $\{i,i+1\}$ are 
 occupied by $X$-particles  and one of them is
$X_{\lbla}$. 
The clocks  $\{N^{i\to i\pm 1}\}$
cannot interchange the $X$-particles across the 
 edge $\{i,i+1\}$ because in the $(\xi,\zeta)$-graphical construction
these are lower priority $\zeta$-particles that do not jump on
top of each other.  The otherwise missing jumps are now
supplied by  the ``new'' clocks that govern the jumps 
described in rules (i)--(iii). 

Combining (i)--(iii)  we can read that
if $X_{\lbla}=i+1$ and  $X_{\lbla-1}=i$, 
then $\lbla$ jumps to $\lbla-1$ with rate $p$. 
This is the first case of type 3 jumps above,
 corresponding to a $\zeta^-$-particle moving from
$i$ to $i+1$ with rate $p$, and the second class particle
$X_{\lbla}$ yielding. 
On the other hand,  if  $X_{\lbla}=i$ and  $X_{\lbla+1}=i+1$ then 
$\lbla$ jumps to $\lbla+1$ with rate $q$. This is
the second case in type 3, 
corresponding to a $\zeta^-$-particle moving from
$i+1$ to $i$ with rate $q$ and exchanging places with
 the second class particle
$X_{\lbla}$. 

We have  verified that the process $(\zeta^-, X_{\lbla})$ 
operates with the correct rates. 

To argue from the rates to the correct distribution of the 
process, we can make use of the process  $(\zeta^-, \zeta)$.
 The processes  $(\zeta^-, X_{\lbla})$ 
and  $(\zeta^-, \zeta)$ determine each other uniquely. 
The virtue of  $(\zeta^-, \zeta)$ is that it has a compact
state space and only nearest-neighbor jumps
with bounded rates.  Hence by the basic  theory
of semigroups and generators of particle systems as developed
in \cite{ips}, 
given the initial configuration,
the distribution of the process  is uniquely determined by the
action of the generator on local functions.  Thus
it suffices to check that individual jumps have the 
correct rates across each edge $\{i,i+1\}$. This  is 
exactly what we did above in the language
of  $(\zeta^-, X_{\lbla})$.
\end{proof}

Similar argument shows that $(\xi, X_{\lblb})$ 
is a $(p,q)$-ASEP with a second class particle.  
To prove 
 Theorem \ref{th:newcoupling} take 
 $Q^\zeta=X_{\lbla}$ and $Q^\xi= X_{\lblb}$. 
This gives the coupling whose existence is claimed in the 
theorem.  

To conclude, let us observe  that the four
processes $(\xi, \xi^+,\zeta^-,\zeta)$ are not in
basic coupling.   For example,  the jump from state
\[
\begin{bmatrix}  \zeta_i&\zeta_{i+1}\\
\zeta^-_i&\zeta^-_{i+1}\\ \xi^+_i&\xi^+_{i+1}\\ 
\xi_i&\xi_{i+1}\end{bmatrix} 
=\begin{bmatrix} 1&1\\0&1\\1&0\\0&0\end{bmatrix}
\quad\text{to state}\quad
\begin{bmatrix} 1&1\\1&0\\0&1\\0&0\end{bmatrix}
\]
happens at rate $q$ (second case of rule (i)), while
in basic coupling this move is impossible.

\section{Biased random walk in an inhomogeneous environment}
\label{sc:RW}
This section proves an estimate that is needed for the
proof of Theorem \ref{th:2classmom}.
Let $Z(t)$ be a continuous-time nearest-neighbor 
random walk on  state-space
$S\subseteq\bZ$ that contains $\bZ_-=\{\dotsc,-2,-1,0\}$. Initially $Z(0)=0$.
 $Z$ attempts to jump from $x$ to $x+1$ with rate $p$ for $x\le -1$,
and from $x$ to $x-1$ with rate $q$ for $x\le 0$.
Assume $p>q=1-p$ and let $\theta=p-q$.  The rates on $S\smallsetminus\bZ_-$ need not
be specified. 
 
Whether jumps are permitted or not is determined by a fixed
environment expressed in terms of $\{0,1\}$-valued 
functions $\{u(x,t): x\in S,\, 0\le t<\infty\}$.  
A jump across edge $\{x-1,x\}$  in either direction is permitted
at time $t$ if $u(x,t)=1$, otherwise not.  
In other words, $u(x,t)$ is the
indicator of the event that edge $\{x-1,x\}$ is open at time $t$. 

\medskip

{\bf Assumption.}  Assume that for all $x\in S$ and  $T<\infty$,
$u(x,t)$ flips between $0$ and $1$ only finitely many times
during $0\le t\le T$.  Assume for convenience
 right-continuity: $u(x,t+)=u(x,t)$.

\begin{lemma} For all $t\ge 0$ and  $k\ge 0$,
\[
P\{ Z(t) \le -k\} \le e^{-2\theta k}.
\]
This bound holds for any fixed environment $\{u(x,t)\}$
subject to the assumption above. 
\label{lm:RW}\end{lemma}

\begin{proof} 
 Let $Y(t)$ be a walk that operates exactly
as $Z(t)$ on $\bZ_-$ but is rectricted to remain in $\bZ_-$
by setting the  rate of jumping from $0$ to $1$ to zero. 
Give $Y(t)$ geometric initial distribution
\[ P\{Y(0)=-j\}= \pi(j)\equiv
 \Bigl(1-\frac{q}{p}\Bigr)\Bigl(\frac{q}{p}\Bigr)^j
\quad\text{for $j\ge 0$.}\]
The initial points satisfy $Y(0)\le Z(0)$ a.s.
Couple the walks through Poisson clocks so that the inequality
 $Y(t)\le Z(t)$ is preserved for all time $0\le t<\infty$.

Without the inhomogeneous  environment $Y(t)$ would be a stationary, 
reversible  birth and death
process.  We argue that even with the environment the time marginals
$Y(t)$ still have distribution $\pi$. 
This suffices for the conclusion, for then
\[
P\{Z(t)\le -k\}\le P\{Y(t)\le -k\}=(q/p)^k
=\exp\bigl(k\log\tfrac{1-\theta}{1+\theta}\bigr)
\le e^{-2\theta k}.
\]

To justify the claim about $Y(t)$, consider approximating processes
$Y^{(m)}(t)$, $m\in\bN$, with the same initial value
$Y^{(m)}(0)=Y(0)$.  $Y^{(m)}(t)$ evolves  so that the environments
$\{u(x,t)\}$ restrict its motion only on edges
$\{x-1,x\}$ for $-m+1\le x\le 0$.  In other words, for walk $Y^{(m)}(t)$
we set $u(x,t)\equiv 1$ for $x\le -m$ and $0\le t<\infty$. 
We couple the walks together so that $Y(t)=Y^{(m)}(t)$ 
until the first time one of the walks exits the interval
$\{-m+1,\dotsc,0\}$.  

Fixing $m$ for a moment,
let $0=s_0<s_1<s_2<s_3<\dotsc$ be a partition of the 
time axis so that $s_j\nearrow\infty$ and the environments
$\{u(x,t): -m<x\le 0\}$ are constant on each interval 
$t\in[s_i,s_{i+1})$.  Then on each time interval $[s_i,s_{i+1})$
$Y^{(m)}(t)$ is a continuous time Markov chain with time-homogeneous
jump rates
\begin{align*}
c(x,x+1)&=\begin{cases} pu(x+1,s_i), &-m\le x\le 0\\
                        p, &x\le -m-1 \end{cases}\\
\intertext{and}
c(x,x-1)&=\begin{cases} qu(x,s_i), &-m+1\le x\le 0\\
                        q, &x\le -m. \end{cases}
\end{align*}
One can check that detailed balance $\pi(x)c(x,x+1)=\pi(x+1)c(x+1,x)$
 holds for all $x\le -1$.  Thus $\pi$ is a reversible measure
for walk $Y^{(m)}(t)$ on each  time interval $[s_i,s_{i+1})$,
and we conclude that 
 $Y^{(m)}(t)$ has distribution $\pi$ for all $0\le t<\infty$. 

The coupling ensures that $Y^{(m)}(t)\to Y(t)$ almost surely
as $m\to\infty$,
and consequently also $Y(t)$  has distribution $\pi$ for all 
$0\le t<\infty$. 
\end{proof}

\section{Moment bounds for the second class particle}
\label{sc:2momproof}
In this section we prove Theorem \ref{th:2classmom}. 

\subsection{Upper bound} 
\label{sec:UB1}
We prove the upper bound
of Theorem \ref{th:2classmom} by proving bounds for 
tail probabilities. We do this first for the right tail
of $Q(t)$. 
Throughout we assume fixed rates  $p> q=1-p$ and abbreviate 
\be
\theta=p-q. 
\label{def:theta}\ee
 Introduce also the notation
\be
\Psi(t)=\Ev^\rho\lvert Q(t)-V^\rho t\rvert.
\label{def:Psi}\ee

\begin{lemma} Let $r\ge 1\vee 8\sqrt\theta$.
Then for each density  $0< \rho< 1$ there exists a  constant 
$C(\rho)\in(0,\infty)$  such that, for
  $t\ge 1$ and $u\ge rt^{2/3}$, 
\be
\Pv^\rho\{Q(t)\ge V^\rho t+u\} \le 
 C(\rho)\Bigl(\,\frac{t^2}{u^4}\Psi(t) + \frac{t^2}{u^3} \Bigr). 
\label{eq:UBge1/2}\ee
The constant $C(\rho)$ is  continuous in $\rho\in(0,1)$,
and $\lim_{\rho\searrow 0}C(\rho)=\infty$. 
\label{lm:UB1}
\end{lemma}

\begin{proof} 
 Assume for convenience that  
$u$ is a positive integer.  Since
$\fl{u}\ge u/2$ for $u\ge 1$, \eqref{eq:UBge1/2}
extends from integers $u$ to real $u$  by an adjustment
of  the constant $C$. 

Fix a density
 $0<\rho<1$ and an auxiliary  density $0<\lambda<\rho$
that will vary in the argument.
Start with the basic coupling of three exclusion processes
$\om\ge\om^-\ge \eta$ with this initial set-up: 
 
(a) 
Initially $\{\om_i(0):i\ne 0\}$ are i.i.d.\ Bernoulli($\rho$) 
distributed 
and $\om_0(0)=1$.

(b) 
Initially $\om^-(0)=\om(0)-\delta_0$.   

(c) Initially
 variables $\{\eta_i(0):i\ne 0\}$ 
are i.i.d.\ Bernoulli($\lambda$) and $\eta_0(0)=0$. 
 The   coupling of the 
initial occupations is 
such that $\om_i(0)\ge\eta_i(0)$ for all $i\ne 0$.

Recall that  basic coupling meant that  these processes
obey  common  Poisson clocks. 

Let $Q(t)$ be the position of the single second class 
particle between $\om(t)$ and $\om^-(t)$, initially
at the origin.  Let $\{X_i(t):i\in\bZ\}$ be the positions
of the $\om-\eta$ second class particles, initially 
labeled so that 
\[
\dotsm<X_{-2}(0)<X_{-1}(0)<X_0(0)=0<X_1(0)<X_2(0)<\dotsm 
\]
These second class particles preserve their labels in the
dynamics and stay ordered.  Thus the $\om(t)$ configuration
consists of first class particles (the $\eta(t)$ process)
and second class particles (the $X_j(t)$'s). 
$\Pv$ denotes the probability measure under which all these
coupled processes live. Note that the marginal distribution
of $(\om,\om^-,Q)$ under $\Pv$ is exactly as it would
be under $\Pv^\rho$. 

For $x\in\bZ$,  
 $J^\om_x(t)$ is  the net current in the 
$\om$-process between space-time positions 
$(1/2,0)$ and $(x+1/2,t)$.  Similarly 
 $J^\eta_x(t)$  in the $\eta$-process, 
and  $J^{\om-\eta}_x(t)$ is  the net current of
second class particles.  
Current 
in the $\om$-process is a sum of  the first class particle
current and the second class particle current:  
\be
J^\om_x(t)= J^\eta_x(t)+J^{\om-\eta}_x(t).
\label{eq:currsum}\ee

 $Q(t)$ is included among the $\{X_j(t)\}$ for all time
because  the basic coupling preserves the
coordinatewise ordering
$\om^-(t)\ge\eta(t)$.   Define the label $\Qlb(t)$ by
$Q(t)=X_{\Qlb(t)}(t)$ with initial value $\Qlb(0)=0$. 
We insert a bound on the label. 

\begin{lemma} 
For all $t\ge 0$ and $k\ge 0$,
\[
\Pv\{\Qlb(t)\ge k\} \le e^{-2\theta k}.
\]
\label{lm:mQ}\end{lemma} 
\begin{proof}[Proof of Lemma \ref{lm:mQ}]
In the basic coupling 
 the label $\Qlb(t)$ evolves as follows.  When $X_{\Qlb-1}$ is
adjacent to $X_{\Qlb}$, $\Qlb$ jumps down by one at rate $p$. 
And when  $X_{\Qlb+1}$ is
adjacent to $X_{\Qlb}$, $\Qlb$ jumps up by one at rate $q$. 
When $X_{\Qlb}$ has no  $X$-particle in either neighboring site,
the label $\Qlb$ cannot jump. 
Thus the situation is like that in  Lemma \ref{lm:RW} 
(with a reversal of lattice directions) 
with  environment  given
by the adjacency of $X$-particles:  
$u(m,t)=\ind\{X_m(t)=X_{m-1}(t)+1\}$. 
However,  the basic coupling mixes together the
evolution of the environment and the walk $\Qlb$,
so the environment is not specified in advance
as required by  Lemma \ref{lm:RW}. 

We can get around this difficulty by imagining an alternative
but distributionally equivalent 
construction for the joint process $(\eta, \om^-,\om)$. 
Let $(\eta,\om)$ obey basic coupling with the given Poisson
clocks $\{N^{x\to x\pm 1}\}$ attached to 
directed edges $(x,x\pm 1)$.  Divide
the $\om-\eta$ particles further into
class II  consisting  of the 
particles   $\om^--\eta$  and 
class III that consists only of the single 
particle  $\om-\om^-=\delta_Q$.  Let class II have 
priority over class III.  
 Introduce another independent set of Poisson
clocks $\{\widetilde N^{x\to x\pm 1}\}$, also attached
to directed edges $(x,x\pm 1)$ of the space $\bZ$ where
particles move.  Let clocks $\{\widetilde N^{x\to x\pm 1}\}$  
 govern the exchanges between 
classes II and III.  In other words, for each edge
$\{x,x+1\}$  clocks
$\widetilde N^{x\to x+ 1}$ and $\widetilde N^{x+1\to x}$ 
are observed if 
sites $\{x,x+1\}$ are both occupied by $\om-\eta$
particles.  All other jumps are prompted by the 
original clocks. 

The rates for individual jumps are the same in this
alternative construction as in the earlier one where
all processes were together in basic coupling.  Thus
the same distribution for the process  $(\eta, \om^-,\om)$
 is created.  

To apply  Lemma \ref{lm:RW} perform the construction
in two steps. First construct the process  $(\eta,\om)$
for all time.  This determines the environment 
$u(m,t)=\ind\{X_m(t)=X_{m-1}(t)+1\}$.  Then run the 
dynamics of classes II and III in this environment.
Now Lemma \ref{lm:RW} gives the bound for $\Qlb$. 
\end{proof}

The preliminaries are ready and we begin 
 to develop a series of inequalities. Let
$u$ and $k$ be positive integers. 
\begin{align}
&\Pv\{Q(t)\ge V^\rho t +u\} \nn\\
&\qquad \le \Pv\{\Qlb(t)\ge k\} 
+ \Pv\{  J^\om_{\fl{V^\rho t}+u}(t) \; - \; J^\eta_{\fl{V^\rho t}+u}(t)
>-k\}. 
\label{line2}\end{align}
To explain the inequality above,
if $Q(t)\ge V^\rho t +u$ and 
$\Qlb(t)<k$ then $X_k(t)>\fl{V^\rho t}+u$. This puts
the bound  \[J^{\om-\eta}_{\fl{V^\rho t}+u}(t)> -k\] on the
second class particle current, because 
at most particles  $X_1,\dotsc, X_{k-1}$ could have
made a negative contribution to this current.

Lemma \ref{lm:mQ} takes care of
the first probability on line \eqref{line2}.
We work on the second probability on line
\eqref{line2}.  

Here is a simple observation that will be used repeatedly.  
Process $\om$ can be coupled
with a stationary density-$\rho$ process $\om^{(\rho)}$ 
 so that the coupled pair $(\om,\om^{(\rho)})$ has 
at most 1 discrepancy.  In this coupling
\be
\lvert J^\om_{x}(t)-J^{\om^{(\rho)}}_{x}(t)\rvert\le 1. 
\label{eq:coupJ}\ee 
This way we can  use computations for 
stationary processes at the expense of small errors.  

 Recall
that $V^\rho=\flux'(\rho)$. Let
$c_1$ below be a constant that absorbs the errors
from using means of stationary processes  
 and from ignoring  integer parts. It satisfies 
$\abs{c_1}\le 3$.  
\begin{align}
\Ev J^\om_{\fl{V^\rho t}+u}(t) - \Ev J^\eta_{\fl{V^\rho t}+u}(t)
&= 
t\flux(\rho)-(\flux'(\rho)t+u)\rho -t\flux(\lambda)+(\flux'(\rho)t+u)\lambda
+c_1\nn\\
&= -\tfrac12 t\flux''(\rho)(\rho-\lambda)^2-u(\rho-\lambda)+c_1 \nn\\
&=t\theta(\rho-\lambda)^2-u(\rho-\lambda)+c_1.  \label{line1.8}
\end{align}
For more general fluxes with nonvanishing $\flux''(\rho)$ the
Taylor expansion would produce more terms above. 

The discussion  splits into three cases according to
the range of $u$-values.  Only the first case requires substantial
work. 

\medskip

{\bf Case 1.} $rt^{2/3}\le u\le \rho\theta t$.

\medskip

Choose
\be
\lambda=\rho-\frac{u}{2\theta t}
\quad\text{and}\quad
k=\left\lfloor\frac{u^2}{8\theta t}\right\rfloor-3.
\label{eq:lakchoice}\ee
The assumptions $r\ge 1\vee 8\sqrt\theta$, $t\ge 1$
 and $u\ge rt^{2/3}$  ensure 
that 
\be
k\ge \frac{u^2}{16\theta t}\ge 1. 
\label{eq:kprop}\ee
In the next line below
 the $-3$ in $k$ absorbs $c_1$ from line \eqref{line1.8}.
Recall that  $\widetilde X=X-EX$ stands for a centered random
variable.  We continue bounding the second probability
 from line \eqref{line2}.
\begin{align}
&\Pv\{  J^\om_{\fl{V^\rho t}+u}(t) \; - \; J^\eta_{\fl{V^\rho t}+u}(t)
>-k\}\nn\\
&\le \Pv\Bigl\{  \widetilde J^\om_{\fl{V^\rho t}+u}(t) \; - \; 
\widetilde J^\eta_{\fl{V^\rho t}+u}(t) \ge  
\frac{u^2}{8t\theta} \Bigr\}\nn\\
&\le \frac{64\theta^2t^2}{u^4} 
\Vv\Bigl\{ J^\om_{\fl{V^\rho t}+u}(t) \; - 
\; J^\eta_{\fl{V^\rho t}+u}(t)\Bigr\}
\nn\\
&\le \frac{128\theta^2t^2}{u^4} 
\Bigl(\Vv\bigl\{ J^\om_{\fl{V^\rho t}+u}(t)\bigr\} 
\; + \; \Vv\bigl\{J^\eta_{\fl{V^\rho t}+u}(t)\bigr\}\,\Bigr).\label{line6}
\end{align} 

We develop bounds on  the variances above, 
first for $J^\om$. Pass to the stationary density-$\rho$ 
 process via \eqref{eq:coupJ} and apply 
\eqref{goal1} of Theorem \ref{th:formulas}:
\begin{align}
\Vv\bigl\{ J^\om_{\fl{V^\rho t}+u}(t)\bigr\}
&\le 2\Var^\rho\bigl\{ J_{\fl{V^\rho t}+u}(t)\bigr\}+2 \nn\\
&=
2\rho(1-\rho) 
\Ev\bigl\lvert Q(t)-\fl{V^\rho t}-u\bigr\rvert + 2\nn\\
&\le
\Ev\lvert Q(t)-V^\rho t\,\rvert + u +3.\nn
\end{align}  
As was already pointed out, as far as $Q(t)$ goes,
the $\Ev^\rho$ expectation in the right member of 
\eqref{goal1} is the same as  $\Ev$ in the coupling
of this section.  Recall  the notation $\Psi(t)$ from
\eqref{def:Psi}, 
bound $3$ by $3u$ (recall that $u\ge 1$),  
and write the above bound in the form 
\be
\Vv\bigl\{ J^\om_{\fl{V^\rho t}+u}(t)\bigr\} \le
\Psi(t) + 4u.
\label{eq:aux7}\ee

In order to get the same bound for 
$\Vv\bigl\{J^\eta_{\fl{V^\rho t}+u}(t)\bigr\}$ we utilize
the coupling developed in Section \ref{sc:newcoup}.
Let $\Var^\lambda$ denote variance in the stationary 
density-$\lambda$  process
and let  $Q^\eta(t)$ denote the position
of a second class particle added to a process 
$\eta$ defined as at the beginning of Section \ref{sec:UB1}. 

Reasoning as was done for \eqref{eq:aux7}:
switch to a stationary density-$\lambda$ process
and  apply  \eqref{goal1}:
\begin{align}
&\Vv\bigl\{ J^\eta_{\fl{V^\rho t}+u}(t)\bigr\}
\le 2\Var^\lambda\bigl\{ J_{\fl{V^\rho t}+u}(t)\bigr\}+2 \nn\\
&\qquad \le 
\Ev^\lambda\bigl\lvert Q^\eta(t)-\fl{V^\rho t}-u\bigr\rvert + 2\nn\\
&\qquad \le
\Ev^\lambda\lvert Q^\eta(t)-V^\rho t\,\rvert + 4u \nn\\
\intertext{Introduce process 
$(\zeta^-(t), Q^\zeta(t), \eta(t), Q^\eta(t))_{t\ge 0}$
coupled as  in Theorem \ref{th:newcoupling}, where 
$\zeta$ starts with Bernoulli($\rho$) occupations away
from the origin and  initially $Q^\zeta(0)=Q^\eta(0)=0$. 
Below apply the triangle inequality and use inequality
   $Q^\zeta(t)\le Q^\eta(t)$ from   Theorem \ref{th:newcoupling}. 
Thus continuing from above: }
&\qquad =
\Ev\lvert Q^\eta(t)-Q^\zeta(t)+Q^\zeta(t)-V^\rho t\,\rvert + 4u\nn\\
&\qquad \le
\Ev\bigl\{ Q^\eta(t)-Q^\zeta(t)\bigr\} +
\Ev\lvert Q^\zeta(t)-V^\rho t\,\rvert + 4u\nn\\
&\qquad= 
V^\lambda t-V^\rho t + \Psi(t) + 4u\nn\\
&\qquad = 2\theta t(\rho-\lambda) + \Psi(t) + 4u\nn\\
&\qquad = \Psi(t) + 5u. 
\label{eq:aux9}
\end{align}  
  Marginally the process
$(\zeta, Q^\zeta)$ is the same as the process
$(\om, Q)$ in the coupling of this section, 
hence the appearance of $\Psi(t)$ above. Then 
we used \eqref{eq:EQH'} for the expectations of the second class
particles and  the  choice 
\eqref{eq:lakchoice} of $\lambda$. 

Insert bounds \eqref{eq:aux7} and \eqref{eq:aux9}
into \eqref{line6} to get
\begin{align}
\Pv\{  J^\om_{\fl{V^\rho t}+u}(t) \; - \; J^\eta_{\fl{V^\rho t}+u}(t)
>-k\}
\le C\theta^2
\Bigl(\,\frac{t^2}{u^4}\Psi(t) + \frac{t^2}{u^3}\,\Bigr).
\label{eq:aux13} \end{align}
where $C=1152$. 
  By \eqref{eq:kprop} and  Lemma \ref{lm:mQ} 
\[\Pv\{\Qlb(t)\ge k\} \le e^{-u^2/8t}\]
which we bound by (recalling $u\ge 1$)
\[
 e^{-u^2/8t}\le C
\biggl(\frac{t}{u^2}\biggr)^2\le C\frac{t^2}{u^3}.
\]
Insert this and \eqref{eq:aux13} into line \eqref{line2} 
to get the upper tail bound
\be 
\Pv\{Q(t)\ge V^\rho t+u\} \le 
 C(\theta)\Bigl(\,\frac{t^2}{u^4}\Psi(t) + \frac{t^2}{u^3}
\Bigr)
\label{eq:aux14}\ee
which is exactly the goal \eqref{eq:UBge1/2}
for {\bf Case 1}. 

\medskip 

{\bf Case 2.}  $\rho\theta t\le u\le 3t$.

\medskip

This comes from the previous case. With $u\le 3t$,
apply 
\eqref{eq:aux14} to $v=\frac13\rho\theta u\le\rho\theta t$
to get the bound
\begin{align}
&\Pv\{Q(t)\ge V^\rho t+u\} \le 
\Pv\{Q(t)\ge V^\rho t+v\} \nn\\
&\qquad \le C(\theta)
\Bigl(\,\frac{t^2}{v^4}\Psi(t) + \frac{t^2}{v^3}
\Bigr)  \le \frac{C(\theta)}{\rho^4}
\Bigl(\,\frac{t^2}{u^4}\Psi(t) + \frac{t^2}{u^3}
\Bigr).
\label{eq:aux15}
\end{align}

\medskip

{\bf Case 3.}  $u\ge 3t$.

\medskip

This comes from a large deviation bound. 
$Q(t)$ is stochastically dominated by a rate 1
Poisson process $Z_t$ and $\abs{V^\rho}\le 1$.
A straightforward exponential Chebyshev 
argument for the Poisson distribution gives the bound
\be
\Pv^\rho\{Q(t)\ge V^\rho t+ u\} 
\le \Pv\{Z_t\ge 2u/3\}\le e^{-Bu} 
\label{eq:aux15.5}\ee
for $t\ge 1$ and  $u\ge 3t$
 with a constant $B$ independent of all the parameters.
  The rightmost member 
of \eqref{eq:aux15.5} is dominated by $Ct^2/u^3$ with constant $C$
a fixed number.  So it can be covered by the rightmost member
of \eqref{eq:aux15} by increasing the constant $C(\theta)$
because $\rho^{-4}\ge 1$.
Renaming the constant in \eqref{eq:aux15}   
gives the conclusion  of Lemma \ref{lm:UB1}. 
\end{proof} 

We extend the bound to both tails. 

\begin{lemma} Let $r\ge 1\vee 8\sqrt\theta$.
Then for each density  $0< \rho< 1$ there exists a  constant 
$C_1(\rho)\in(0,\infty)$  such that, for
  $t\ge 1$ and $u\ge rt^{2/3}$, 
\be
\Pv^\rho\{\,\lvert Q(t)- V^\rho t\rvert \ge u\} \le 
 C_1(\rho)\Bigl(\,\frac{t^2}{u^4}\Psi(t) + \frac{t^2}{u^3}\,\Bigr). 
\label{eq:aux16}\ee
The constant $C_1(\rho)$ is  continuous in $\rho\in(0,1)$,
and $\lim_{\rho\to\{0,1\}}C_1(\rho)=\infty$. 
\label{lm:UB2}
\end{lemma}

\begin{proof} Consider coupled processes $(\om,\om^+)$ 
that start  with i.i.d.~density $\rho\in(0,1)$
occupations away from the
origin and initially $\om^+(0)=\om(0)+\delta_0$.
So there is one second class  particle
$Q$ initially at the origin.  Define new processes
\[
\omhat_i(t)=1-\om_i(t)\quad\text{and}\quad
\omhat^-_i(t)=1-\om^+_i(t).
\]
Process 
$\omhat$ records the dynamics
of holes in the $\om$-process, and similarly for $\omhat^-$. 
The discrepancy between $(\omhat^-,\omhat)$ is exactly 
the same as the discrepancy between $(\om,\om^+)$. 
That is,  $\Qhat(t)=Q(t)$. 
Processes $(\omhat^-,\omhat)$ are instances of
$(\phat,\qhat)$-ASEP with density $\rhohat=1-\rho$
and rates $\phat=q$ and $\qhat=p$. 

To recover the original rates $(p,q)$ 
 we  reflect the lattice across
the origin.  Define 
\[
\om^R_i(t)=\omhat_{-i}(t)\quad\text{and}\quad
\om^{R-}_i(t)=\omhat^-_{-i}(t).
\]
Process $(\om^{R-},\om^R)$ is an instance of 
$(p,q)$-ASEP at density $1-\rho$ with a discrepancy, 
so the previous bound
\eqref{eq:aux15} applies. The discrepancy is now
at $Q^R(t)=-\Qhat(t)=-Q(t)$. The characteristic speed is
$V^{1-\rho}=-V^\rho$.   Hence also 
\[
\Psi(t)=\Ev^\rho\lvert Q(t)-V^\rho t\rvert
=\Ev^\rho\lvert Q^R(t)-V^{1-\rho} t\rvert.
\]
By an application of  inequality \eqref{eq:aux15} to the
process $(\om^{R-},\om^R, Q^R)$, 
\begin{align}
&\Pv^\rho\{Q(t)\le V^\rho t-u\} 
=\Pv^{\rho}\{ Q^R(t)\ge V^{1-\rho} t+u\}\nn\\
&\qquad \le \frac{C(\theta)}{(1-\rho)^4}
\Bigl(\,\frac{t^2}{u^4}\Psi(t) + \frac{t^2}{u^3}\,\Bigr).
\label{eq:aux15.1}
\end{align}
 Combining 
\eqref{eq:aux15} and \eqref{eq:aux15.1} we have the 
conclusion \eqref{eq:aux16} with constant 
\be
C_1(\rho)=\frac{C(\theta)}{\rho^4\wedge(1-\rho)^4}
\label{eq:UBconstants}\ee
where $C(\theta)$ is a constant that depends on $\theta$.  
\end{proof}

Let us also record the two-sided large deviation bound,
one side of which was argued above in \eqref{eq:aux15.5}.

\begin{lemma}  There exists a constant $B$ such that 
\be
\Pv^\rho\{\,\lvert Q(t)- V^\rho t\rvert \ge  u\} 
\le e^{-Bu} 
\label{eq:aux17}\ee
for $t\ge 1$,  $u\ge 3t$  and $0<\rho<1$.
\label{lm:ld1}\end{lemma} 

Inequalities \eqref{eq:aux16} and \eqref{eq:aux17} 
 give the upper bound on the moments
of the second class particle in \eqref{eq:mainubd}
via  a two-step integration argument.  We keep track of the
precise constants for a while for use in the lower bound
proof to come. 
First by \eqref{eq:aux16} 
\begin{align*}
\Psi(t)&= \int_0^\infty 
\Pv^\rho\{\,\lvert Q(t) -V^\rho t\rvert \ge u\}\,du\\
&\le 
rt^{2/3}
+ C_1(\rho) \int_{rt^{2/3}}^\infty
\Bigl(\,\frac{t^2}{u^4}\Psi(t) + \frac{t^2}{u^3}\,\Bigr)\,du  \\
&\le \frac{C_1(\rho)}{3r^{3}}\Psi(t) + 
\Bigl(r+\frac{C_1(\rho)}{2r^{2}} \Bigr) t^{2/3}.
\end{align*}
 Fixing $r= \max\{1,\, 8\sqrt\theta,\, C_1(\rho)^{1/3}\}$ 
shows that $\Psi(t)\le C_2(\rho)t^{2/3}$ for a new
constant $C_2(\rho)$ with the same properties as $C_1(\rho)$. 

Put this back into the estimate \eqref{eq:aux16}
 to get
\be 
\Pv^\rho\{\,\lvert Q(t) -V^\rho t\rvert \ge u\} \le 
 C_3(\rho)\Bigl(\,\frac{t^{8/3}}{u^4} + \frac{t^2}{u^3}\Bigr)
\le 
 C_3(\rho)\frac{t^2}{u^3}
\label{eq:aux18}\ee
with new constant $C_3(\rho)$. The second inequality used 
$u\ge rt^{2/3}$. 
Now take $\mom>1$ 
and this time use both  \eqref{eq:aux16} and \eqref{eq:aux17}: 
\begin{align}
&\Ev^\rho\lvert Q(t)-V^\rho t\rvert^\mom 
=\mom \int_0^\infty \Pv^\rho\{\,\lvert Q(t) -V^\rho t\rvert \ge u\}
u^{\mom-1}\,du\nn\\
&\le 
r^\mom t^{2\mom/3} 
\;+\; C_3(\rho) \int_{rt^{2/3}}^{3t}
{t^2}{u^{\mom-4}}\,du
\;+\;\mom\int_{3t}^\infty e^{-Bu}u^{\mom-1}\,du.
\label{eq:aux18.7} 
\end{align}
Performing the integrations gives these bounds:
\[
\Ev^\rho\lvert Q(t)-V^\rho t\rvert^\mom  \le 
\begin{cases} 
\dfrac{C_4(\rho)}{3-\mom} t^{2\mom/3}  
&\text{for $1<\mom<3$ and $t\ge 1$,}\\[8pt]
C_5(\rho)t^2\log t  &\text{for $\mom=3$  and $t\ge e$,}\\
C_6(\rho,\mom) t^{\mom-1}  &\text{for $3<\mom<\infty$ and $t\ge 1$.}\\
\end{cases} \]
One can see that  the constants $C_i$  are continuous in
 $\rho$.  
$C_6$ diverges to $\infty$ as $\mom\searrow 3$ 
or $\mom\nearrow\infty$.
   
After renaming constants, this gives 
 the upper bounds in Theorem
\ref{th:2classmom} and \eqref{eq:addub}. 
We record one case here for use in the lower bound proof. 

\begin{proposition} 
 There exists a continuous function 
$0<C_{UB}(\rho)<\infty$ of density 
$0<\rho<1$ such that 
for $t\ge 1$, 
\[ \Ev^\rho\lvert Q(t)-V^\rho t\rvert \le C_{UB}(\rho)t^{2/3}. 
\]
\label{pr:UB1}\end{proposition}

\subsection{Lower bound}
Fix a density $0<\rho<1$.
For the lower bound we prove that $t^{-2/3}\Psi(t)$ has
a positive lower bound for large enough $t$ where
$\Psi(t)$ was defined in \eqref{def:Psi}. 
The lower bound \eqref{eq:mainlbd} follows then for all $m\ge 1$
from Jensen's inequality. 

Let $a_1, a_2$ be finite positive constants and set 
\be
b=a_2^2/(32\theta).
\label{eq:bfroma_2}\ee
In the argument $a_1$ and $a_2$ will be chosen
sufficiently large relative to
$p$ and relative to the constants 
$C_{UB}(\cdot)$ in Proposition \ref{pr:UB1}.  
Define an auxiliary density 
$\lambda=\rho-bt^{-1/3}$. 
Fix $t_0\ge 1$ and  large enough so that 
$\lambda\in(\rho/2,\rho)$ for $t\ge t_0$. 
 Restrict all subsequent calculations
to $t\ge t_0$.  
Define  positive
integers 
\be
u=\fl{a_1t^{2/3}} \quad\text{and}\quad 
n=\fl{V^\lambda t}-\fl{V^\rho t}+u.
\label{def:n}\ee

Start with a basic coupling of
three processes $\eta\le\eta^+\le \zeta$:  

\medskip

(a) Initially $\eta$ has
i.i.d.~Bernoul\-li($\lambda$) occupations 
$\{\eta_i(0):i\ne -n\}$ and $\eta_{-n}(0)=0$. 

(b)  Initially $\eta^+(0)=\eta(0)+\delta_{-n}$. 
$Q^{(-n)}(t)$ is the location of the unique
discrepancy between $\eta(t)$ and $\eta^+(t)$. 

(c)  Initially $\zeta$ has independent 
occupation variables, coupled with $\eta(0)$
as follows: 

(c.1) $\zeta_i(0)=\eta_i(0)$ for $-n<i\le 0$.

 (c.2) $\zeta_{-n}(0)=1$.

 (c.3)   For  $i< -n$ and $i>0$ variables $\zeta_i(0)$ 
are i.i.d.~Bernoulli($\rho$) and 
$\zeta_i(0)\ge\eta_i(0)$. 

 Thus the initial density of $\zeta$ is piecewise
constant:  on the segment 
$\{-n+1,\dotsc,0\}$ $\zeta(0)$ is i.i.d.\ with
density $\lambda$, at site $-n$ $\zeta(0)$ has density $1$,
and elsewhere on $\bZ$ $\zeta(0)$ is i.i.d.\ with
density $\rho$.

\medskip

Label the $\zeta-\eta$ second class particles
as $\{Y_m(t):m\in\bZ\}$ 
so that initially 
\[
\dotsm< Y_{-1}(0)<Y_{0}(0)= -n=Q^{(-n)}(0)<0 < Y_{1}(0)<Y_{2}(0)<\dotsm
\]
Let again $\Qlb(t)$ be the label such that 
 $Q^{(-n)}(t)=Y_{\Qlb(t)}(t)$. 
Initially $\Qlb(0)=0$. 
The inclusion 
$Q^{(-n)}(t)\in\{Y_m(t)\}$ 
  persists for all
time because the basic coupling preserves the 
ordering $\zeta(t)\ge\eta^+(t)$. 
Through the basic coupling
 $\Qlb$ jumps to the left
with rate $q$ and to the right with rate $p$, but only
when there is a $Y$-particle  adjacent to $Y_{\Qlb}$.
As in the proof of Lemma \ref{lm:mQ} we can apply
Lemma \ref{lm:RW} to prove this statement: 

\begin{lemma}
For all $t\ge 0$ and $k\ge 0$,
\[
\Pv\{\Qlb(t)\le -k\} \le e^{-2\theta k}.
\]
\label{lm:mQ1}\end{lemma}

Fix $a_1>0$ large enough so that 
$u=\fl{a_1t^{2/3}} $ 
satisfies, for all the auxiliary densities  $\lambda$, 
\be
\Pv\{ Q^{(-n)}(t)\ge \fl{V^\rho t}\}
=\Pv\{ Q^{(-n)}(t)\ge -n+\fl{V^\lambda t} +u\} 
 \le \tfrac12.
\label{eq:ucond1}\ee
This is possible because the constant $C_{UB}(\lambda)$ in
the  upper bound in
Proposition \ref{pr:UB1} is continuous in the density and   
 we restricted to $t\ge t_0$ to ensure that 
$\lambda\in(\rho/2,\rho)$.  

If $Q^{(-n)}(t)\le \fl{V^\rho t}$ and  $\Qlb(t)>-k$ then 
 $Y_{-k}(t)\le  \fl{V^\rho t}$. As was argued for \eqref{line2}
in Section \ref{sec:UB1} 
this  implies a bound on the second class particle current:
\[
J^\zeta_{\fl{V^\rho t}}(t)-J^\eta_{\fl{V^\rho t}}(t)
=J^{\zeta-\eta}_{\fl{V^\rho t}}(t)\le k.
\]
 So from \eqref{eq:ucond1} 
\be\begin{split}
\tfrac12 &\le \Pv\{ Q^{(-n)}(t)\le \fl{V^\rho t}\}\\
&\le \Pv\{\Qlb(t)\le -k\}
+\Pv\{ J^\zeta_{\fl{V^\rho t}}(t)-J^\eta_{\fl{V^\rho t}}(t) \le k\}. 
\end{split}\label{eq:aux36}\ee

Take $a_2>0$ large enough so that 
\be a_2\ge 8+8\sqrt{C_{UB}(\alpha)} 
\quad\text{for $\alpha\in(\rho/2,\rho)$.}
\label{eq:a2kappa}\ee
Increase $t_0$ further   so that for $t\ge t_0$ 
\be
\Pv\{\Qlb(t)\le -a_2t^{1/3}+3\} \le \tfrac14
\label{eq:aux38}\ee
(utilizing Lemma \ref{lm:mQ1}). 
Combine displays \eqref{eq:aux36} and \eqref{eq:aux38} 
 with $k=\fl{a_2t^{1/3}}-2$
to get  the next inequality. Then split the probability. 
\begin{align}
\tfrac14&\le 
\Pv\{ J^\zeta_{\fl{V^\rho t}}(t)-J^\eta_{\fl{V^\rho t}}(t) \le 
a_2t^{1/3}-2\} \nn\\
&\le 
\Pv\{ J^\zeta_{\fl{V^\rho t}}(t) \le 2a_2t^{1/3} 
+t\theta(2\rho\lambda-\lambda^2) \} 
\nn\\
&\qquad + 
\Pv\{ J^\eta_{\fl{V^\rho t}}(t) \ge  a_2t^{1/3} 
+t\theta(2\rho\lambda-\lambda^2)+2\}. 
\label{line19}
\end{align}

We treat line \eqref{line19}. Recall that  $P^\lambda$ 
denotes probabilities of  a stationary 
density-$\lambda$ process. As described above \eqref{eq:coupJ}
we can imagine  a basic coupling  in which the 
$\eta$-process differs from a  stationary 
density-$\lambda$ process
by at most one discrepancy. Compute  the mean current in
the stationary process: 
\begin{align*}
E^\lambda\{ J_{\fl{V^\rho t}}(t)\bigr\} 
&= tH(\lambda)-\lambda\fl{V^\rho t} \\
&\le tH(\lambda)-\lambda{V^\rho t}+1
=t\theta(2\rho\lambda-\lambda^2)+1. 
\end{align*}
From these comes the bound
\begin{align}
\text{line \eqref{line19}} 
&\le 
P^\lambda\{ J_{\fl{V^\rho t}}(t) \ge  a_2t^{1/3} 
+t\theta(2\rho\lambda-\lambda^2)+1\}\nn\\
&\le P^\lambda\{ \widetilde J_{\fl{V^\rho t}}(t) \ge  a_2t^{1/3} \}
\le a_2^{-2}t^{-2/3}
\Var^\lambda\bigl\{ J_{\fl{V^\rho t}}(t)\bigr\}\nn\\
&\le \frac{\Ev^\lambda\lvert Q(t)-\fl{V^\rho t}\rvert}{a_2^{2}t^{2/3}}
\le \frac{\Ev^\lambda\lvert Q(t)-V^\lambda t\rvert}{a_2^{2}t^{2/3}}
+ \frac{2\theta bt^{2/3}+1}{a_2^{2}t^{2/3}}\nn\\
&\le C_{UB}(\lambda)a_2^{-2} + \tfrac1{16} + \tfrac1{64}\le \tfrac18.
\label{line22}
\end{align} 
After Chebyshev above  we used 
 \eqref{goal1} and introduced a second class particle 
$Q(t)$ in a density-$\lambda$ system under the measure
 $\Pv^\lambda$. 
 $\fl{V^\rho t}$ was replaced with $V^\lambda t$ at the cost
of an error $1$ for dropping integer parts, and 
 $V^\rho -V^\lambda =2\theta bt^{-1/3}$. Last
we applied  the upper
bound from   Proposition
\ref{pr:UB1}   and 
properties \eqref{eq:bfroma_2},  \eqref{eq:a2kappa} and $t\ge t_0\ge 1$. 

Put this last bound back into line \eqref{line19}. This leaves
\be
\tfrac18\le 
\Pv\{ J^\zeta_{\fl{V^\rho t}}(t) \le 2a_2t^{1/3} 
+t\theta(2\rho\lambda-\lambda^2) \}
\label{line24}\ee

To treat this probability 
 we take the estimation back
to a stationary density-$\rho$ process by inserting the
Radon-Nikodym factor.  Let $\gamma$
denote the distribution of the initial $\zeta(0)$ 
configuration described by (a)--(c) in the beginning
of this section. As before $\nu^\rho$ is the
density-$\rho$ i.i.d.~Bernoulli measure.  Their Radon-Nikodym
derivative is 
\[
f(\om)=\frac{d\gamma}{d\nu^\rho}(\om)= 
\frac1{\rho}\ind\{\om_{-n}=1\}\cdot
\prod_{i=-n+1}^0\Bigl( \frac{\lambda}{\rho}\ind\{\om_i=1\}
+ \frac{1-\lambda}{1-\rho}\ind\{\om_i=0\}\Bigr).
\]
Recalling that $n=O(t^{2/3})$ and $\rho-\lambda=O(t^{-1/3})$, 
\be
E^\rho(f^2)=\frac1\rho\Bigl(1+\frac{(\rho-\lambda)^2}{\rho(1-\rho)}
\Bigr)^n 
\le \rho^{-1}e^{n(\rho-\lambda)^2/\rho(1-\rho)}
\le c_1(\rho)^2,
\label{def:c1}\ee
where $c_1(\rho)$ is 
 independent of $t$, continuous in $\rho$
 but diverges to $\infty$ as $\rho\to\{0,1\}$. 

Let $A$ denote the exclusion process event
\[
A=\{ J_{\fl{V^\rho t}}(t) \le 2a_2t^{1/3} 
+t\theta(2\rho\lambda-\lambda^2) \}.
\]
  Then from \eqref{line24}
\begin{align}
\tfrac18&\le \Pv\{\zeta\in A\}
= \int P^\om(A)\,\gamma(d\om) 
= \int P^\om(A)f(\om)\,\nu^\rho(d\om)\nn\\
&\le \bigl(P^\rho(A)\bigr)^{1/2}\bigl(E^\rho(f^2)\bigr)^{1/2}
\le c_1(\rho) \bigl(P^\rho(A)\bigr)^{1/2}.
\label{line26} 
\end{align}

In the final calculation
 we bound $P^\rho(A)$ by Chebyshev to return to the 
current variance. Note that 
\begin{align*}
E^\rho\bigl\{J_{\fl{V^\rho t}}(t)\bigr\} 
= t\flux(\rho)-\rho\fl{V^\rho t}
=t\theta\rho^2+\rho V^\rho t-\rho\fl{V^\rho t}
\ge t\theta\rho^2.
\end{align*}
Also, if we choose $a_2$ large enough, namely $a_2>(2048\theta)^{1/3}$,
we can ensure that 
\be
c_2\equiv b^2\theta-2a_2>0.
\label{eq:ba2}\ee 
From line \eqref{line26} we have, utilizing  
$\lambda=\rho-bt^{-1/3}$ and \eqref{eq:ba2}, 
\begin{align*}
(8c_1(\rho))^{-2}&\le P^\rho(A)= P^\rho\{ J_{\fl{V^\rho t}}(t) \le 2a_2t^{1/3} 
+t\theta(2\rho\lambda-\lambda^2) \}\\
&\le P^\rho\{ \widetilde J_{\fl{V^\rho t}}(t) \le 2a_2t^{1/3} 
-t\theta(\rho-\lambda)^2 \}\\
&\le P^\rho\{ \widetilde J_{\fl{V^\rho t}}(t) 
\le -(b^2\theta-2a_2)t^{1/3} \}\\
&\le c_2^{-2}t^{-2/3}
\Var^\rho \bigl\{J_{\fl{V^\rho t}}(t)\bigr\}\\
&\le c_2^{-2}t^{-2/3}\Psi(t)
\end{align*}
by \eqref{goal1} again and the abbreviation
\eqref{def:Psi}. 
We have  assumed $t\ge t_0$.
The first and last
lines of the calculation above show that 
$t^{-2/3}\Psi(t)$ has a positive lower bound for 
$t\ge t_0$.  
This lower bound depends on $\rho$ through $c_1(\rho)$ 
in \eqref{def:c1}, and vanishes as $\rho\to\{0,1\}$.

The lower bound of Theorem \ref{th:2classmom} is proved. 
\bibliography{refsmarton}
\bibliographystyle{plain}
\end{document}